\documentclass[10pt]{amsart}
\usepackage{amssymb}
\usepackage{amsthm}

\renewcommand{\phi}{\varphi}
\newcommand{\C}{\mathbb{C}}

\newcommand{\Kth}{K_\theta}
\newcommand{\Sth}{S_\theta}
\newcommand{\T}{\mathbb T}

\newcommand{\Z}{\mathbb Z}
\newcommand{\D}{\mathbb D}
\newcommand{\Tau}{\mathcal{T}_\theta}
\renewcommand{\le}{\leqslant}
\renewcommand{\ge}{\geqslant}
\newcommand{\ov}{\overline}
\newcommand{\der}{{\rm D}}
\newcommand{\DER}{\bar{{\rm D}}}
\newcommand{\lder}{{\rm \partial}}
\newcommand{\LDER}{{\rm \bar{\partial}}}
\newcommand{\ko}{k_0}
\newcommand{\kot}{\tilde{k}_0}
\newcommand{\kl}{k_\lambda}
\newcommand{\klt}{\tilde{k}_\lambda}
\newcommand{\km}{k_\mu}
\newcommand{\kmt}{\tilde{k}_\mu}
\newcommand{\kernel}{\mathop{\mathrm{Ker}}\nolimits}
\newcommand{\ran}{\mathop{\mathrm{Ran}}\nolimits}

\newcommand{\spn}{\mathop{\mathrm{span}}\nolimits}

\newtheorem{Thm}{Theorem}

\newtheorem{Lem}{Lemma}[section]
\newtheorem{Prop}[Lem]{Proposition}

\newtheorem{Rem}[Lem]{Remark}

\newtheorem*{Prob}{Problem}
\def\beginpf{\medskip\noindent {\bf Proof.} \;}

\title[Truncated Toeplitz operators of finite rank]
{Truncated Toeplitz operators of finite rank}

\dedicatory{Accepted for publication in Proceedings of the\,AMS}

\author{R.~V.~Bessonov}

\address{St.Petersburg Department of Steklov Mathematical Institute RAS, Fontanka 27, St.Petersburg, 191023, Russia}

\email{bessonov@pdmi.ras.ru}

\thanks{
This work is partially supported by the RFBR grant 11-01-00584-a, by the V.A.Rokhlin grant 2012 and by the Chebyshev Laboratory  (Department of Mathematics and Mechanics, St. Petersburg State University) under RF Government grant 11.G34.31.0026
}

\begin{document}

\begin{abstract}
We give a complete description of the finite-rank truncated Toeplitz operators.
\end{abstract}

\maketitle
%
%
%
%
\section{Introduction}\label{intro}
Truncated Toeplitz operators are compressions of multiplication operators on\,\,$L^2$ to model subspaces $\Kth = H^2\ominus\theta H^2$ of the Hardy class $H^2$, where $\theta$ is an inner function. In \cite{sar} D.Sarason gave a characterization of all truncated Toeplitz operators of rank one. Moreover, he constructed a class of finite-rank truncated Toeplitz operators and asked whether this class exhausts all truncated Toeplitz operators having finite rank. Our main result, Theorem \ref{main}, answers this question in the affirmative.

\subsection{Definitions.} As usual, we identify the Hardy class $H^2$ in the unit disk $\D$ with the subspace of the space $L^2$ on the unit circle\,$\T$ via non-tangential boundary values. A function $\theta\in H^2$ is called inner if $|\theta(z)| =1$ almost everywhere with respect to the Lebesgue measure on $\T$. Denote by $P_\theta$ the orthogonal projection from $L^2$ onto $\Kth$. The truncated Toeplitz operator $A_\phi$ with symbol $\phi \in L^2$ is the mapping
\begin{equation}\label{eq-1}
A_\phi: f \mapsto P_{\theta} (\phi f), \quad f \in \Kth\cap L^\infty.
\end{equation}
We deal only with bounded truncated Toeplitz operators. For $\lambda \in \D$ define
\begin{equation}\label{eq-26}
\kl = \frac{1-\ov{\theta(\lambda)}\theta}{1 - \ov{\lambda}z}, \qquad \klt = \frac{\,\theta - \theta(\lambda)}{z-\lambda}.
\end{equation}
The function $\kl$ is the reproducing kernel of the space $\Kth$ at the point $\lambda$; $\klt$ is the conjugate kernel at $\lambda$ (see Section \ref{shift} for details). For an integer $n \ge 0$, denote by $\Omega_n = \Omega(\theta, n)$ the set of all points $\lambda \in \T$ such that every function from $\Kth$ and its derivatives up to order $n$ have non-tangential limits at $\lambda$. A description of $\Omega_n$ is given by P.R.Ahern and D.N. Clark in \cite{ac70}, we discuss it in Section\;\ref{shift}. In particular, for $\lambda \in \T$, we have $\lambda \in \Omega_0$ if and only if $\kl, \klt \in \Kth$.

\medskip

The only compact Toeplitz operator on the Hardy space $H^2$ is the zero operator\,\cite{halmos}. The situation is different for Toeplitz operators on $\Kth$. \pagebreak It was proved in \cite{sar} that the general rank-one truncated Toeplitz operator on $\Kth$ is a scalar multiple of
\begin{equation}\label{eq-2}
\kl \otimes \klt \qquad \mbox{or} \qquad \klt \otimes \kl,
\end{equation}
for some $\lambda \in \D \cup \Omega_0$, where we use the standard notation for rank-one operators in the Hilbert space: $x \otimes y: h \mapsto (h,y)x$. Some finite-rank truncated Toeplitz operators can be obtained from the rank-one operators \eqref{eq-2} by differentiation. Consider the analytic mapping $\Phi: \lambda \mapsto \klt\otimes\kl$ from the unit disk $\D$ into the space of all bounded operators on $\Kth$. Take an integer $n \ge 1$. Denote by $\der^n [\klt\otimes\kl]$ the value $\Phi^{(n)}(\lambda)$ of the $n$-th derivative of $\Phi$ at the point $\lambda \in \D$. For $\lambda \in \Omega_n$ let $\der^n [\klt\otimes\kl]$ denote the $n$-th angular derivative of $\Phi$ at $\lambda$. Existence of this derivative follows from results of P.R.Ahern and D.N.Clark and from representation \eqref{eq-31} of $\der^n [\klt\otimes\kl]$ in terms of derivatives of reproducing kernels; see Section \ref{kernels} for details. Since the set of all truncated Toeplitz operators is linear space closed in the weak operator topology \cite{sar}, the derivative $\der^n [\klt\otimes\kl]$ determines truncated Toeplitz operator of rank $n$ for every $\lambda \in \D\cup\Omega_n$. Parallel arguments work for the adjoint operator $\DER^n [\kl\otimes\klt]$: it is the derivative of $\kl\otimes\klt$ of order $n$ with respect to $\bar\lambda$. Operators of the form
\begin{equation}\label{eq-4}
\DER^n[\kl\otimes\klt], \qquad \der^n [\klt\otimes\kl], \qquad n\ge 0, \;\; \lambda \in \D \cup \Omega_n,
\end{equation}
were originally constructed in \cite{sar} as examples of finite-rank truncated Toeplitz operators that are not linear combinations of the rank-one operators \eqref{eq-2}.
\subsection{The main result.}\label{statement}
\begin{Thm}\label{main}
The general finite-rank truncated Toeplitz operator on $\Kth$ is a finite linear combination of the operators in \eqref{eq-4}.
\end{Thm}
The key step in proving of Theorem \ref{main} is the identification of range of a general finite-rank truncated Toeplitz operator. Set
$$
F(\lambda,n)=\begin{cases}
\ran\DER^n[\kl\otimes\klt],&\text{if $\lambda \in \D\cup\Omega_n$;}\\
\ran\der^n[\tilde{k}_{\lambda^{*}}\otimes k_{\lambda^{*}}], &\text{if $\lambda \in \D_e$,}
\end{cases}
$$
 where $\D_e = \{z: \;|z| > 1\} \cup\{\infty\}$, $\lambda^* = 1/\bar{\lambda}$ and $\infty^* =0$. It will be shown in Section\,\ref{mteta} that a subspace $F\subset\Kth$ is range of an operator in \eqref{eq-4} if and only if we have $F=F(\lambda,n)$ for some $n\ge 0$ and $\lambda \in\D\cup\Omega_n\cup\D_e$. Theorem\,\ref{main} will be proved as soon as we establish the following three results.
\begin{Lem}\label{lem-range} Suppose $A$ is a finite-rank truncated Toeplitz operator on $\Kth$. Then there exists a finite collection of points $\lambda_k \in \D\cup\Omega_{n_k}\cup\D_e$ such that
\begin{equation}\label{eq-44}
\ran A = F(\lambda_1, n_1) \dotplus \ldots \dotplus F(\lambda_s, n_s).
\end{equation}
\end{Lem}
\begin{Lem}\label{lem-splitting}
Suppose $A$ is a truncated Toeplitz operator on $\Kth$ with range of the form \eqref{eq-44}. Then $A$ is a sum of $s$ truncated Toeplitz operators $A_k$ such that $\ran A_k = F(\lambda_k, n_k)$.
\end{Lem}
\begin{Lem}\label{lem-elementary}
Suppose $A$ is a truncated Toeplitz operator on $\Kth$ with the range $\ran A = F(\lambda,n)$, where $\lambda \in \D\cup\Omega_n\cup\D_e$. Then $A$ is a finite linear combination of the operators in \eqref{eq-4}.
\end{Lem}

In Section \ref{shift} we collect some standard results about the space $\Kth$. Section \ref{range} concerns a description of range of a general bounded truncated Toeplitz operator.  Lemmas \ref{lem-range}, \ref{lem-splitting}, \ref{lem-elementary} will be proved in Section \ref{proof}.

\subsection{A duality approach.}
The rank-one truncated Toeplitz operators can be regarded as a point evaluation in a special Banach space of analytic functions. Given an inner function $\theta$, define the space $X_a$ by
$$
X_a = \Bigl\{\sum_{0}^{\infty} x_k y_k: \; x_k,y_k \in \Kth \;\mbox{ and }\; \sum_{0}^{\infty} \|x_k\|\!\cdot\! \|y_k\| < \infty\Bigr\}.
$$
This space was introduced in \cite{bbk}, where the fact that $X_a$ is the predual of the space of all bounded truncated Toeplitz operators established. The pairing is given by
\begin{equation}\label{eq-43}
\bigl\langle A, \sum x_k y_k\bigr\rangle = \sum (A x_k, y_k).
\end{equation}
It was shown in \cite{bbk}, that $\langle \klt\otimes\kl, xy\rangle= (x\tilde{y})(\lambda)$, where $\tilde{y} = \bar{z}\theta \bar{y}\in \Kth$. This yields the nice formula
$$\bigl\langle \der^n [\klt\otimes\kl], \sum x_k y_k\bigr\rangle = \Bigl(\sum x_k \tilde{y}_k\Bigr)^{(n)}(\lambda).$$
A similar relation holds for the operators $\DER^n [\kl\otimes\klt]$.  In the sense of the pairing  \eqref{eq-43} the operator $\kl\otimes\klt$ is the point evaluation at  $\lambda^* = 1/\bar{\lambda}$ of functions from $X_a$ up to the scalar $(\lambda^*/\theta(\lambda^*))^2$, and $\DER^n [\kl\otimes\klt]$ is its derivative with respect to $\bar{\lambda}$. Since $|\lambda^*|\ge 1$, the point evaluation should be understood in terms of pseudocontinuations of functions from $\Kth$ to the exterior of the unit disk (see Lecture\,II in \cite{nik}).
\subsection{Notations.}
\begin{itemize}
\item[] $\Z_+$ is the set of nonnegative integers;
\item[] $\Tau$ is the linear space of all bounded truncated Toeplitz operators on $\Kth$;
\item[] $\kernel A$ is kernel of a bounded operator $A$;
\item[] $\ran A$ is range of a bounded operator $A$;
\item[]$\ov{\ran A}$ is the closure of $\ran A$;
\item[] $\langle f_1 \ldots f_n \rangle=\spn\{f_j,\; j=1\ldots n\}$;
\item[] $H^{\bot}$ is the orthogonal complement to a subspace $H$;
\item[] $H_1 + H_2$ is the linear span of the union $H_1 \cup H_2$;
\item[] $H_1 \dotplus H_2$ is the direct sum of two subspaces $H_1$, $H_2$.
\end{itemize}

\section{Preliminaries}\label{shift}
\noindent This section contains the preliminary information concerning spaces $\Kth$ and truncated Toeplitz operators: reproducing kernels, Clark unitary perturbations, Sarason's characterization of truncated Toeplitz operators.  A more detailed discussion is available in \cite{sar}, \cite{nik},  \cite{cimmat}.
\subsection{The conjugation.}\label{conjugation} The space $\Kth$ is closed under the conjugation
\begin{equation}\label{eq-45}
C: x \mapsto \bar{z}\theta\bar{x}.
\end{equation}
Truncated Toeplitz operators are complex symmetric with respect to $C$, which means $CA = A^*C$, see \cite{sar}. Hence the ranges of $A$ and $A^*$ are mutually conjugate for every  operator $A \in \Tau$: $\ran A^* = C\ran A$. For more information on this property of truncated Toeplitz operators see \cite{gar4}, \cite{gar5}.
\subsection{Reproducing kernels and their derivatives.}\label{kernels} For $\lambda \in \D$, set
\begin{equation}\label{eq-5}
\kl=\frac{1 - \ov{\theta(\lambda)}\theta}{1-\bar{\lambda}z}, \qquad
\klt=\frac{\,\theta - \theta(\lambda)}{z - \lambda}.
\end{equation}
Note that $\klt = C\kl$, where the conjugation $C$ is defined by \eqref{eq-45}. The function $\kl$ (respectively, $\klt$) is the reproducing kernel (respectively, conjugate reproducing kernel) of the space $\Kth$ at the point\,$\lambda$:
\begin{equation}\label{eq-6}
(f, \kl) = f(\lambda), \qquad (f, \klt) = \ov{(Cf)(\lambda)}.
\end{equation}
Differentiating \eqref{eq-6}, we obtain
\begin{equation}\label{eq-7}
(f, \LDER^n\kl) = f^{(n)}(\lambda), \qquad (f, \lder^n\klt) = \ov{(Cf)^{(n)}(\lambda)}.
\end{equation}
In what follows the symbols $\LDER^n\kl$ and $\lder^n\klt$ denote the $n$-th derivatives of $\kl, \klt$ with respect to $\bar\lambda, \lambda$, respectively. For example, in the case $n=1$ we have
$$
\LDER\kl = \lim_{\mu \to \lambda}\frac{\kl - \km}{\ov{\lambda} - \ov{\mu}},
\qquad
\lder\klt = \lim_{\mu \to \lambda} \frac{\klt - \kmt}{\lambda - \mu}.
$$
Let $\theta = B_{\Lambda}S_{\nu}$ be an inner function with the set of zeroes $\Lambda = (a_k)_{k=1}^{N}$, repeated according to multiplicity, and the singular part $S_\nu$ that corresponds to a singular measure $\nu$ on the unit circle $\T$. Take a point $\lambda \in \T$ and an integer $n \in \Z_+$. The following result is in \cite{ac70}, see also \cite{cl}.
\begin{Thm}[P.R.Ahern and D.N.Clark]\label{ac-theorem}
The following are equivalent:
\begin{itemize}
\item[(a)] functions $f$, $f', \ldots f^{(n)}$ have non-tangential limits at $\lambda$ for every $f \in \Kth$;
\item[(b)] the non-tangential limits $\LDER^j\kl = \lim_{\mu \to \lambda} \LDER^j\km$ and $\lder^j\klt = \lim_{\mu \to \lambda} \lder^j \kmt$ exist in norm of $\Kth$ for every $j = 0 \ldots n$;
\item[(c)] $\sum_{k=1}^{N} \frac{1 - |a_k|^2}{|1 - \lambda \bar a_k|^{2n+2}} < \infty$ and $\int_{\T} \frac{d\nu(\xi)}{|1 - \lambda \bar\xi|^{2n+2}} < \infty$.
\end{itemize}
\end{Thm}
Denote by $\Omega_n =\Omega(\theta, n)$ the set of all points $\lambda \in \T$ that satisfy the conditions of Theorem \ref{ac-theorem}. For $\lambda \in \D$ we have
\begin{equation}\label{eq-31}
\begin{split}
&A_\phi = \DER^n[\kl \otimes \klt] = \sum_{k=0}^{n}C_{n}^{k}\left(\LDER^k\kl\otimes\lder^{n-k}\klt\right), \qquad
\phi(z)=\frac{n!\cdot\ov{\theta(z)}}{(\bar z-\bar \lambda)^{n+1}};\\
&A_\phi = \der^n[\klt \otimes \kl] = \sum_{k=0}^{n}C_{n}^{k}\left(\lder^k\klt\otimes\LDER^{n-k}\kl\right), \qquad
\phi(z)=\frac{n!\cdot\theta(z)}{(z-\lambda)^{n+1}}.
\end{split}
\end{equation}
by using the Leibniz formula for derivatives of a bilinear expression. It follows from Theorem \ref{ac-theorem} that formulas \eqref{eq-5}-\eqref{eq-31} hold in the sense of non-tangential boundary values for every $\lambda \in \Omega_n$. Thus, the operators in \eqref{eq-4} are exactly the operators in\,\eqref{eq-31} for $n\in\Z_+$ and $\lambda \in \D \cup \Omega_n$.
\subsection{The restricted shift}\label{mteta} Consider the operator $\Sth: f \mapsto P_\theta (zf)$. We have
\begin{equation}\label{eq-11}
\begin{split}
\Sth\kl = (\kl -\ko)/\bar\lambda; \qquad &\Sth\klt = \lambda\klt -\theta(\lambda)\ko;\\
\Sth\ko = \LDER\ko; \qquad &\Sth\kot = -\theta(0)\ko,
\end{split}
\end{equation}
where $\lambda \neq 0$ is a point from $\D \cup \Omega_0$. In particular, if $\theta(0) = 0$ then $\kot \in \kernel \Sth$. Actually, we have $\kernel\Sth = \langle \kot \rangle$ in that case.
\begin{Prop}\label{lem-nder}
We have $\Sth^n\ko=\frac{1}{n!}\LDER^n\ko$ and  $\Sth\lder^n\kot = n\lder^{n-1}\kot - \theta^{(n)}(0)\ko$ for every integer $n \ge 1$.
\end{Prop}
\beginpf Take a function $f \in \Kth$ and consider $(f, \Sth^n\ko) = ((\Sth^*)^nf, \ko)$. The space $\Kth$ is invariant under the backward shift operator $S^*: f \mapsto (f - f(0))/z$, see \cite{nik}. Hence $S^*|\Kth = \Sth^*$ and $((\Sth^*)^nf, \ko) = \frac{1}{n!}f^{(n)}(0)$. It follows from \eqref{eq-7} that $f^{(n)}(0) = (f, \LDER^n\ko)$, and therefore $(f, \Sth^n\ko) = \frac{1}{n!}(f, \LDER^n\ko)$. Since this equality holds for every function $f \in \Kth$, we obtain $\Sth^n\ko=\frac{1}{n!}\LDER^n\ko$. Differentiating the identity $\Sth\klt = \lambda\klt -\theta(\lambda)\ko$ with respect to $\lambda$ at the point $\lambda=0$, we get the formula $\Sth\lder^n\kot = n\lder^{n-1}\kot - \theta^{(n)}(0)\ko$. \qed

\medskip

\noindent Consider the subspaces $F(\lambda,n)$ defined in Section \ref{statement}. It follows from formula \eqref{eq-31} that
\begin{equation}\label{eq-36}
F(\lambda,n)=\begin{cases}
\spn\{\LDER^j\kl,\; j=0\ldots n\},&\text{if $\lambda \in \D\cup\Omega_n$;}\\
\spn\{\lder^j \tilde{k}_{\lambda^{*}},\; j=0\ldots n\}, &\text{if $\lambda \in \D_e$.}
\end{cases}
\end{equation}
For each point $\lambda \in \Omega_0$ we have $\lambda^*=\lambda$ and $\klt = \bar{\lambda}\theta(\lambda)\kl$. Therefore, 	
$$\spn\{\LDER^j\kl,\; j=0\ldots n\} = \spn\{\lder^j \tilde{k}_{\lambda^{*}},\; j=0\ldots n\}$$
for all $\lambda \in \Omega_n$ and $n \in \Z_+$. We now see from \eqref{eq-31} that range of the operator ${\DER^n[\kl \otimes \klt]}$ coincides with range of the operator $\der^n[\klt \otimes \kl]$ if $\lambda \in \Omega_n$. Therefore, a subspace $F\subset\Kth$ is range of an operator in \eqref{eq-4} if and only if we have $F=F(\lambda,n)$ for some $n\ge 0$ and $\lambda \in\D\cup\Omega_n\cup\D_e$, as claimed in Section\,\ref{statement}.
\begin{Prop}\label{lem-1}
We have
\begin{equation}\label{eq-16}
\begin{split}
&\Sth F(\lambda,n) \subset F(\lambda,n) \dotplus \langle\ko\rangle,\mbox{ if } \lambda \neq 0;\\
&\Sth F(0,n) \subset F(0,n) \dotplus \langle\LDER^{n+1}\ko\rangle.
\end{split}
\end{equation}
\end{Prop}
\beginpf The first formula in \eqref{eq-16} can be obtained from \eqref{eq-11} by differentiation. The second one follows from Proposition \ref{lem-nder}. \qed

\subsection{The Frostman shift.}
Let $\theta$  be an inner function. The Frostman shift of $\theta$ corresponding to the point $\theta(0)$ is the inner function $\Theta = \frac{\theta - \theta(0)}{1- \ov{\theta(0)}\theta}$. We have $\Theta(0)=0$.
Define the unitary operator $J: \Kth \to K_\Theta$ by
$$
J: f\mapsto \frac{\sqrt{1-|\theta(0)|^2}}{1- \ov{\theta(0)}\theta}f.
$$
\begin{Prop}\label{prp12}
For $\lambda \in \D$ and $n \in \Z_+$ we have
\begin{equation}\label{eq-29}
J\spn\{\LDER^j\kl,\; j=0\ldots n\} = \spn\{\LDER^j\kl^\Theta,\; j=0\ldots n\},
\end{equation}
where $\kl^\Theta$ is the reproducing kernel of the space $K_\Theta$ at the point $\lambda$.
\end{Prop}
\beginpf The formula $J \kl = \bigl(1 - \theta(0)\ov{\theta(\lambda)}\bigr)/\bigl(\sqrt{1-|\theta(0)|^2}\bigr) \kl^\Theta$ follows from the definition of $J$ and implies \eqref{eq-29} by  differentiation; see details in Section 13 of \cite{sar}.\qed

\medskip

\noindent The following fact is a particular case of Theorem 13.2 from \cite{sar}.
\begin{Prop}\label{prp13}
A bounded operator $A$ is a truncated Toeplitz operator on $\Kth$ if and only if the operator $JAJ^{-1}$ is a truncated Toeplitz operator on $K_\Theta$.
\end{Prop}

\subsection{Clark's unitary perturbations.}\label{clark} In \cite{cl}
 D.N.Clark described one-dimensional unitary perturbations of $\Sth$. Given a number $\alpha \in \T$, define
\begin{equation}\label{eq-8}
U_\alpha = \Sth + c_\alpha \ko \otimes \kot, \qquad c_\alpha = \frac{\alpha + \theta(0)}{1 - |\theta(0)|^2},
\end{equation}
where $\ko,\kot$ are the reproducing kernels \eqref{eq-5} at the origin. The operators $U_\alpha$ are unitary and cyclic; every one-dimensional unitary perturbation of $\Sth$ is $U_\alpha$ for an appropriate number $\alpha \in \T$, see \cite{cl}. It is shown in \cite{cl} that the spectral measure $\sigma_\alpha$ of the unitary operator $U_\alpha$ can be chosen so that
\begin{equation}\label{eq-9}
U_\alpha = V_{\alpha}^{-1}M_zV_\alpha,
\end{equation}
where $V_\alpha: \Kth \to L^2(\sigma_\alpha)$ is the unitary operator that sends functions from a dense subset of $\Kth$ to their boundary values on $\T$; $M_z$ is the operator of multiplication by $z$ on $L^2(\sigma_\alpha)$. A.G.Poltoratski \cite{polt} established the existence of non-tangential boundary values $\sigma_\alpha$-almost everywhere for all functions from $\Kth$. Thus, the operator $V_\alpha$ is the well-defined unitary embedding $\Kth \to L^2(\sigma_\alpha)$. For every function $f \in \Kth$ we have
\begin{equation}\label{eq-32}
f(z) = \int_{\T}(V_\alpha f)(\xi)\frac{1 - \bar\alpha \theta(z)}{1 - \bar\xi z} \, d\sigma_{\alpha}(\xi), \quad z \in \D.
\end{equation}
The following fact is due to D.N.Clark \cite{cl}, see also \cite{cimmat}.
\begin{Prop}[D.N.Clark]\label{prp5} We have $\sigma_\alpha \{\xi \in \T:\; \theta(\xi)=\alpha\} = \sigma_\alpha(\T)$ for each $\alpha \in \T$. The measure $\sigma_\alpha$ has an atom at a point $\lambda \in \T$ if and only if $\lambda \in \Omega_0$ and $\theta(\lambda) = \alpha$. In that case we have $V_{\alpha}^{-1} \mathbb{I}_{\{\lambda\}} = \sigma_\alpha(\{\lambda\})\cdot\kl$, where $\mathbb{I}_{\{\lambda\}}$ denotes the indicator of the singleton $\{\lambda\}$.
\end{Prop}
\begin{Prop}\label{prop27}
If $(z-\lambda)^{-n-1} \in L^2(\sigma_\alpha)$ for some $\lambda \in \T$ and $n \in \Z_+$, then $\lambda \in \Omega_n$ and $\theta(\lambda) \neq \alpha$.
\end{Prop}
\beginpf At first, let $\theta(0) = 0$. In this case constants lie in the space $\Kth$. Formula \eqref{eq-32} with $f \equiv 1$ gives us
\begin{equation}\label{eq-33}
\frac{1}{1 - \bar\alpha \theta(z)} = \int_{\T}\frac{d\sigma_{\alpha}(\xi)}{1 - \bar\xi z}, \quad z \in \D.
\end{equation}
We have
\begin{equation}\label{eq-40}
\frac{1}{|1 - \bar\xi z|} \le \frac{c}{|1 - \bar\xi\lambda|}, \quad \xi \in \T,
\end{equation}
for some constant $c$, as $z$ tends non-tangentially to $\lambda$. Since $(1 - \bar\xi \lambda)^{-1} \in L^2(\sigma_\alpha)$, we see from \eqref{eq-33} and \eqref{eq-40} that the function $\theta$ has non-tangential limit at $\lambda$  and $\theta(\lambda)\neq\alpha$. Similarly, differentiating \eqref{eq-33} with respect to $z$, one can  prove that $\theta'$, $\theta'', \ldots, \theta^{(n)}$ also have non-tangential limits at $\lambda$. In the case $\theta(0) \neq 0$ this fact follows from the consideration of the Frostman shift of the function $\theta$.

Take $f \in \Kth$ and consider $f^{(j)}$, where $0 \le j \le n$ is integer. It follows from \eqref{eq-32} that
\begin{equation}\label{eq-34}
f^{(j)}(z) = \int_{\T}(V_\alpha f)(\xi)\left(\frac{1 - \bar\alpha \theta(z)}{1 - \bar\xi z}\right)^{(j)}  d\sigma_{\alpha}(\xi), \quad z \in \D.
\end{equation}
Since $(1 - \bar\lambda z)^{-j-1} \in L^2(\sigma_\alpha)$, we see from \eqref{eq-40} and \eqref{eq-34} that $f^{(j)}$ has the non-tangential limit at the point $\lambda$ for every $j=0\ldots n$. Thus, we have $\lambda \in \Omega_n$. \qed

\medskip

\noindent We now describe boundary values of functions from subspaces $F(\lambda,n)$ in \eqref{eq-36}.
\begin{Prop}\label{prp7}
Let $\alpha \in \T$, $n \in \Z_+$ and let $\lambda \in \D \cup\Omega_n\cup\D_e$, where in the case $|\lambda| = 1$ we assume that $\theta(\lambda) \neq \alpha$. We have
\begin{gather}
V_\alpha F(\lambda,n) = \spn\{(z - \lambda^*)^{-j}, \; j=1\ldots n+1\}, \quad \lambda \neq 0, \; \lambda^{*} = 1/\bar{\lambda};\label{eq-101}\\
V_\alpha F(0,n) = \spn\{z^j, \; j=0\ldots n\}.\label{eq-102}
\end{gather}
\end{Prop}
\beginpf In the case $|\lambda| \neq 1$ formula \eqref{eq-101} follows from the definition of $V_\alpha$. Formula \eqref{eq-102} is a consequence of Proposition \ref{lem-nder}. Take a point $\lambda = \lambda^*$ from $\Omega_0$ such that $\theta(\lambda) \neq \alpha$. By Proposition \ref{prp5} we have $\sigma_\alpha(\{\lambda\}) = 0$ and
$$
V_\alpha\kl = \frac{1 - \ov{\theta(\lambda)}\alpha}{1 - \bar\lambda z}
$$
$\sigma_\alpha$--almost everywhere.
Since $V_\alpha\kl \in L^2(\sigma_\alpha)$, we have $(z-\lambda)^{-1} \in L^2(\sigma_\alpha)$ and therefore \eqref{eq-101} holds in the case $n=0$. For $n=1$ and $\lambda \in \Omega_1$, consider
$$
V_\alpha\LDER\kl = \frac{-\ov{\theta'(\lambda)}\alpha}{1 - \bar\lambda z} + \frac{z(1 - \ov{\theta(\lambda)}\alpha)}{(1 - \bar\lambda z)^2}.
$$
Since $V\LDER\kl \in L^2(\sigma_\alpha)$  and $(z-\lambda)^{-1} \in L^2(\sigma_\alpha)$, we have $(z-\lambda)^{-2} \in L^2(\sigma_\alpha)$. Hence formula \eqref{eq-101} holds in the case $n=1$. Arguing as above, we prove \eqref{eq-101} for all $\lambda \in \Omega_n$ and $n \in \Z_+$. \qed
\begin{Prop}\label{prp14}
Suppose that an inner function $\theta$ is not a finite Blashke product. Then every finite collection of the functions $\LDER^{s_k}k_{\lambda_{k}}$,  $\lder^{t_{k}}\tilde{k}_{\mu_{k}}$, where $\lambda_k \in \D \cup \Omega_{s_k}$ and $\mu_k \in \D \cup \Omega_{t_k}$, is linearly independent in $\Kth$.
\end{Prop}
\beginpf Since $\theta$ is not a finite Blashke product, the space $\Kth$ has infinite dimension, see Lecture II in \cite{nik}. Hence the space $L^2(\sigma_\alpha)$, $|\alpha|=1$, has infinite dimension as well. The result now follows from  Proposition \ref{prp7}.\qed
\subsection{A characterization.}\label{characterization}
In what follows we will often use the following characterization of truncated Toeplitz operators.
\begin{Thm}[D.Sarason, \cite{sar}]\label{ds}
A bounded operator $A$ on $\Kth$ is a truncated Toeplitz operator if and only if there exist functions $\psi, \chi$ in $\Kth$ such that
\begin{equation}\label{s}
A - \Sth A \Sth^* = \psi\otimes k_0 + k_0\otimes \chi,
\end{equation}
in which case $A = A_{\psi + \bar{\chi}}$.
\end{Thm}
\begin{Rem}\label{rem1} Suppose $\theta(0) = 0$; then $\psi  = A\ko -\ov{\chi(0)}\ko$.
\end{Rem}
\beginpf Apply both sides of \eqref{s} to the vector $\ko\equiv 1$ and use the relation ${\Sth^*1 = 0}$.\qed
\section{The range of a bounded truncated Toeplitz operator}\label{range}
\noindent In this section we prove the following result.
\begin{Prop} \label{prp1} Let $A \in \Tau$, and assume that $\ov{\ran A} \neq \Kth$. Then
\begin{equation} \label{eq1}
\Sth \ov{\ran A}\subset \ov{\ran A} \dotplus \langle \bar{\partial}^n k_0\rangle,
\end{equation}
where $n \in \Z_+$ is the maximal integer such that  $\bar{\partial}^j k_0 \in \ov{\ran A}$ for every $0 \le j < n$.
\end{Prop}

\noindent For the proof we need two lemmas.
\begin{Lem}\label{lem1} Let $\theta(0) = 0$, $A \in \Tau$. Then
\begin{itemize}
\item[(1)] $\Sth A\left[\langle\tilde{k}_0\rangle^\bot\right] \subset \overline{\ran A} + \langle k_0\rangle$;
\item[(2)] $\Sth\ov{\ran A} \subset \overline{\ran A} + \langle k_0, g \rangle$ for some $g \in \Kth$.
\end{itemize}
\end{Lem}
\beginpf Indeed, by Theorem \ref{ds} and Remark \ref{rem1} we have $\Sth A \Sth^* h \in \ov{\ran A} + \langle k_0\rangle$ for every vector $\Sth^* h$ from $\ov{\ran\Sth^*} = \langle \tilde{k}_0\rangle^\bot$. This proves assertion $(1)$. To prove $(2)$, one can take $g = \Sth A\tilde{k}_0$. \qed
\begin{Lem} \label{lem2} Let $\theta(0) = 0$, $A \in \Tau$. Assume that $k_0 \notin \overline{\ran A}$. Then   $$\Sth \ov{\ran A} \subset \ov{\ran A} \dotplus \langle k_0\rangle.$$
\end{Lem}
\beginpf Denote by $P_1, P_2$ the orthogonal projections on $\Kth$ with the ranges $\ov{\ran A^*}$ and $\langle \tilde{k}_0 \rangle$, respectively. Note that $A = AP_1$. For every $h \in \Kth$ we have
$$\Sth Ah =  \Sth AP_1 h = \Sth AP_2P_1 h + f_1,$$
where $f_1 = \Sth AP_2^{\bot}P_1h$. It follows from assertion $(1)$ of Lemma \ref{lem1} that $f_1 \in\overline{\ran A} \dotplus \langle k_0\rangle$. Proceeding inductively, we obtain the relation
\begin{equation}\label{eq-51}
\Sth Ah = \Sth A(P_2P_1)^n h + f_n
\end{equation}
for some $f_n \in \overline{\ran A} \dotplus \langle k_0\rangle$. Since $\ran A^*  = C\ran A$ (see Section \ref{conjugation}) and $k_0 \notin \overline{\ran A}$, we have $\tilde{k}_0 \notin \ov{\ran A^*}$, and therefore $\|P_2P_1\| < 1$. Passing to the limit in \eqref{eq-51}, we see that $\Sth Ah \in \overline{\ran A} \dotplus \langle k_0\rangle$ and the result follows. \qed

\medskip\noindent {\bf Proof of Proposition \ref{prp1}.} \;At first, let $\theta(0) = 0$. In the case $k_0 \notin \overline{\ran A}$, Lemma \ref{lem2} gives us formula \eqref{eq1} with $n = 0$. Now suppose that $k_0 \in \overline{\ran A}$. Since the family $\{\LDER^j k_0\}_{0}^{\infty}$ is complete in $\Kth$ (see formula \eqref{eq-7}), one can choose the maximal integer $n\ge 1$ such that $\LDER^j k_0 \in \ov{\ran A}$ for every $0 \le j < n$ and $\LDER^n k_0 \notin \ov{\ran A}$. It follows from Proposition \ref{lem-nder} and Lemma \ref{lem1} that
\begin{equation}\label{eq3}
\bar{\partial}^n k_0 = n\Sth \bar{\partial}^{n-1}k_0 \in \Sth\ov{\ran A} \subset \ov{\ran A} + \langle k_0, g \rangle = \ov{\ran A} + \langle g \rangle
\end{equation}
for some $g \in \Kth$ such that $\Sth\ov{\ran A} \subset \overline{\ran A} + \langle k_0, g \rangle$. By the construction we have $\bar{\partial}^n k_0 \notin \overline{\ran A}$.	Comparing this with \eqref{eq3}, we get $g \in \ov{\ran A} \dotplus \langle\LDER^n k_0\rangle$ and thus
$$\Sth\ov{\ran A}  \subset \overline{\ran A} + \langle k_0, g \rangle \subset \ov{\ran A} \dotplus \langle\LDER^n k_0\rangle.$$
Hence the Proposition is proved in the case $\theta(0)=0$. The general situation  can be reduced to this case by using Propositions \ref{prp12} and \ref{prp13}.\qed

\medskip

\noindent\textbf{Remark.} Truncated Toeplitz operators are symmetric with respect to the conjugation $C$, see Section \ref{conjugation}. Using this fact, one can obtain an analogue of Proposition \ref{prp1} for $\ov{\ran A}$ with $\Sth^*$ and $\lder^m\kot$ in place of $\Sth$ and $\LDER^n\ko$.

\medskip

Let $U_\alpha = \Sth + c_\alpha\ko\otimes\kot$ be the Clark unitary perturbation of $\Sth$. Consider the embedding $V_\alpha :\Kth \to L^2(\sigma_\alpha)$ from Section \ref{clark}.
\begin{Prop}\label{prp3}
Let $A \in \Tau$, and assume that $\ov{\ran A} \neq \Kth$. Set $F = V_\alpha\ov{\ran{A}}$. Then
$zF \subset F \dotplus \langle z^n\rangle
$, where $n \in \Z_+$ is the maximal integer such that  $z^j \in F$ for every $0 \le j < n$.
\end{Prop}
\beginpf  By Proposition \ref{prp1} we have $\Sth \ov{\ran A}\subset \ov{\ran A} \dotplus \langle \bar{\partial}^n k_0\rangle$, where $n \in \Z_+$ is the maximal integer such that  $\bar{\partial}^j k_0 \in \ov{\ran A}$ for every $0 \le j < n$. Hence $U_\alpha\ov{\ran A} \subset \ov{\ran A} \dotplus \langle\LDER^n\ko\rangle
$. It remains to apply the operator $V_\alpha$ to both sides of this inclusion and use Proposition \ref{prp7}.\qed

\begin{Prob}
Given a finite Borel measure $\nu$ supported on the unit circle $\T$, to describe all subspaces $F \subset L^2(\T, \nu)$ such that $zF \subset F \dotplus \langle 1\rangle$.
\end{Prob}
\noindent In the next section we treat the finite-dimensional case of this problem.

\section{Proof of Theorem \ref{main}}\label{proof}
Hereinafter we assume that the space $\Kth$ has infinite dimension (equivalently, the inner function $\theta$ is not a finite Blashke product). The case ${\rm dim}\,\Kth < \infty$ is considered in \cite{sar}, where the following fact is proved: every truncated Toeplitz operator on $\Kth$, ${\rm dim}\,\Kth < \infty$,  is a finite linear combination of the rank-one operators \eqref{eq-2}.

\medskip

Theorem \ref{main} follows from Lemmas \ref{lem-range}, \ref{lem-splitting}, \ref{lem-elementary}. We now turn to proving this results.
\subsection{Proof of Lemma \ref{lem-range}}\label{fr-range} The proof is based on the following proposition.
\begin{Prop}\label{prp41}
Let $\sigma$ be a finite Borel measure supported on the unit circle $\T$. Suppose that $F \subset L^2(\sigma)$ is a finite-dimensional subspace satisfying $zF \subset F \dotplus \langle 1\rangle$. If $F$ does not contain indicators of singletons, then there exists a finite collection of points $\lambda_k \in \C$ such that $F = Q(\lambda_1, p_1) \dotplus Q(\lambda_2, p_2) \dotplus \ldots \dotplus Q(\lambda_s, p_s)$, where $Q(\lambda_k,p_k) = \spn\{(z-\lambda_k)^{-j},\; j=1 \ldots p_k\}$, $k=1\ldots s$.
\end{Prop}

\beginpf Denote by $\mathcal{P}$ the non-orthogonal projection on $F \dotplus \langle 1\rangle$ with the range $F$ and the kernel $\langle 1\rangle$. Let $M_z$ be the operator of  multiplication by the independent variable on $L^2(\sigma)$. The finite-rank operator $T = \mathcal{P}M_z: F \to F$ has complete family of root vectors. Consider the root subspace $G_\lambda =\kernel(T - \lambda I)^p$, $G_\lambda \neq \kernel(T - \lambda I)^{p-1}$. The Proposition will be proved as soon as we show that $G_\lambda = Q(\lambda,p)$.

\medskip

{\it Case 1.} $\lambda \in \T$, $\sigma(\{\lambda\}) > 0$.

\medskip

\noindent Take a vector $f \in \kernel(T - \lambda I) \subset G_\lambda$. We have $(z-\lambda)f = c $ for a constant $c \in \C$. Therefore, $f$ is a scalar multiple of the indicator of the singleton $\{\lambda\}$. But that contradicts the hypothesis. Hence this case does not arise.

\medskip

{\it Case 2.} $\lambda \in \C$, $\sigma(\{\lambda\}) = 0$.

\medskip

\noindent Take a vector $f \in G_\lambda$ such that  $f_{1} = (T - \lambda I)^{p-1}f \neq 0$. We have $(T - \lambda I)f_{1} = 0$ by the construction. On the other hand, we have $(T - \lambda I)f_{1} = (z-\lambda)f_{1} - c_{1}$ for some constant $c_1 \in \C$. Taking into account $\sigma(\{\lambda\}) = 0$, we see that
$$f_{1} = \frac{c_{1}}{z-\lambda}.$$
If $p=1$, we stop the procedure. Otherwise put $f_{2} = (T - \lambda I)^{p-2}f$ and consider $(T - \lambda I)f_{2} = (z-\lambda)f_{2} - c_{2}$. Since $(T - \lambda I)f_{2} = f_{1}$, we have
$$f_{2} = \frac{c_{1}}{(z-\lambda)^2} + \frac{c_{2}}{z-\lambda}.$$
Continuing this procedure, we obtain
$$f = f_p = \frac{c_{1}}{(z-\lambda)^p} + \frac{c_{2}}{(z-\lambda)^{p-1}} + \ldots + \frac{c_{p}}{z-\lambda}.$$
Since $c_1 \neq 0$ and $f_j \in G_\lambda$ for every $j = 1\ldots p$, we get $(z-\lambda)^{-j} \in G_\lambda$ and hence $Q(\lambda,p) \subset G_\lambda$. Now take an arbitrary vector $f \in G_\lambda$ and find a number $r$, $1 \le r \le p$, such that $f \in \kernel(T - \lambda I)^r$ but $f_{1} = (T - \lambda I)^{r-1}f \neq 0$. Arguing as above, we see that $f \in Q(\lambda,r)\subset  Q(\lambda,p)$ and thus $G_\lambda \subset  Q(\lambda,p)$. \qed

\medskip

\noindent\textbf{Proof (Lemma \ref{lem-range}).} Let $A$ be a finite-rank truncated Toeplitz operator on $\Kth$. Since the subspace $\ran A$ has finite dimension, it cannot contain an infinite system of reproducing kernels, see Proposition \ref{prp14}. Therefore, by Proposition \ref{prp5} we can choose the Clark measure $\sigma_\alpha$ so that the subspace $F = V_\alpha \ran A$ of the space $L^2(\sigma_\alpha)$ does not contain indicators of singletons. It follows from Proposition \ref{prp3} that $zF \subset F \dotplus \langle z^n\rangle$,
where $n \in \Z_+$ is the maximal integer such that  $z^j \in F$ for every $0 \le j < n$. The subspace $\bar{z}^{n}F$ satisfies the assumptions of Proposition \ref{prp41}. We have
\begin{equation*}
F = z^n(Q(\lambda_1, p_1) \dotplus Q(\lambda_2, p_2) \dotplus \ldots \dotplus Q(\lambda_s, p_s)).
\end{equation*}
In the case $n=0$, the application of Propositions \ref{prop27} and  \ref{prp7} concludes the proof. Now  assume that $n\ge 1$. Since $z^j \in F$ for every integer $0 \le j <n$, we necessarily have $\lambda_{k_0}=0$, $p_{k_0}=n$ for some $1\le k_0 \le s$. Renumber the sequence $\{\lambda_k\}$ so that $\lambda_{1} = 0$. The subspace $z^n Q(0, n)$ is the set of polynomials of degree at most $n-1$. A simple algebra gives us
\begin{equation*}
F = z^n Q(0, n) \dotplus Q(\lambda_2, p_2) \dotplus \ldots \dotplus Q(\lambda_s, p_s).
\end{equation*}
Now the result follows from Propositions \ref{prop27} and  \ref{prp7}.  \qed

\subsection{Proof of Lemma \ref{lem-splitting}.} \label{splitting-section}
Before the proof we need the following general result.
\begin{Prop}\label{splitting}
Let $n\in \Z_+$, $A \in \Tau$. Suppose that $\ov{\ran A} = F_1 \dotplus F_2$, where $F_1$ and $F_2$ are subspaces of $\Kth$ such that
\begin{equation} \label{eq-38}
\Sth F_1 \subset F_1 \dotplus \langle\ko\rangle; \qquad
\Sth F_2 \subset F_2 \dotplus \langle \bar{\partial}^n k_0\rangle, \quad \bar{\partial}^j k_0 \in F_2, \quad 0 \le j < n.
\end{equation}
Let also $\LDER^{n}\ko\notin \ov{\ran A}$. Then $A = A_1 + A_2$, where $A_{k} \in \Tau$ and $\ov{\ran A_{k}} = F_{k}$, $k=1,2$.
\end{Prop}
\beginpf Denote by $\mathcal{P}$ the non-orthogonal projection on $F_1 \dotplus F_2 \dotplus \langle\LDER^{n}\ko\rangle$ with the the range $F_2 \dotplus \langle\LDER^{n}\ko\rangle$ and the the kernel $F_1$. We want to show that $A_2 = \mathcal{P}A \in \Tau$. By Theorem\,\ref{ds}, we need to check that
\begin{equation}\label{eq-17}
A_2 - \Sth A_2 \Sth^* = \psi\otimes k_0 + k_0\otimes \chi
\end{equation}
for some $\psi,\chi \in \Kth$. We have
\begin{equation*}
A_2 - \Sth A_2 \Sth^* = \mathcal{P}(A - \Sth A \Sth^*) + (\mathcal{P} \Sth - \Sth \mathcal{P})A \Sth^*.
\end{equation*}
Since $A$ is a truncated Toeplitz operator, it satisfies \eqref{s} with $\psi_1,\chi_1 \in \Kth$. Hence,
$$
\mathcal{P}(A - \Sth A \Sth^*) = (\mathcal{P}\psi_1)\otimes k_0 + (\mathcal{P} k_0)\otimes \chi_1 = (\mathcal{P}\psi_1)\otimes k_0 + k_0\otimes \chi_1.
$$
Next, the operator $\mathcal{P}\Sth - \Sth \mathcal{P}$ vanishes on $F_2$ and maps $F_1$ to $\langle\ko\rangle$. Therefore, we have $\ran(\mathcal{P} \Sth - \Sth \mathcal{P})A \Sth^* \subset \langle\ko\rangle$. This proves\,\eqref{eq-17}. Now put $A_1 = A - A_2$ and obtain the required representation. \qed

\medskip

\noindent\textbf{Proof (Lemma \ref{lem-splitting}).} It follows from Proposition \ref{lem-1} that any splitting of the sum in \eqref{eq-44} into two summands gives us subspaces $F_1, F_2$ with property \eqref{eq-38}. Consequently applying Proposition \ref{splitting}, we obtain the required.

\subsection{Proof of Lemma \ref{lem-elementary}.} \label{elementary}
Let $A$ be a truncated Toeplitz operator on $\Kth$ with the range $\ran A = F(\lambda,n)$, where $\lambda \in \D\cup\Omega_n\cup\D_e$.
Truncated Toeplitz operators are complex symmetric with respect to the conjugation $C$, see Section \ref{conjugation}. Hence
\begin{equation*}
\begin{split}
&\ran A = \spn\{\lder^j\tilde{k}_{\lambda^*},\; j=0\ldots n\} \quad \mbox{in the case }\lambda \in \D_e,\\
&\ran A^{*} = \spn\{\lder^j\klt,\; j=0\ldots n\} \quad \mbox{in the case }\lambda\in\D\cup\Omega_n.
\end{split}
\end{equation*}
Passing if necessary to the adjoint operator, we can assume  that
$$\ran A = \spn\{\lder^j\kmt,\; j=0\ldots n\}$$
for some point $\mu \in \D\cup\Omega_n$. Then $\ran A^* = \spn\{\LDER^j\km,\; j=0\ldots n\}$. Every such operator has the form
\begin{equation}\label{eq-19}
A = \sum_{0 \le p,q\le n} a_{p,q}\bigl(\lder^p\kmt\otimes\LDER^q\km\bigr)
\end{equation}
for some coefficients $a_{p,q} \in \C$. We claim that
$A = \sum_{s =0}^{n} a_{0,s} \cdot \der^s [\kmt\otimes\km]$.
Consider firstly the case $\mu \neq 0$. Set $T_{pq} = \lder^p\kmt\otimes\LDER^q\km$. For $1 \le p,q \le n$ we have
 \begin{equation}\label{eq-20}
 \begin{split}
 T_{pq} - \Sth T_{pq}\Sth^*
 &= \lder^p\kmt\otimes\LDER^q\km - \lder^p(\Sth\kmt)\otimes\LDER^q(\Sth\km)\\
 &= \lder^p\kmt\otimes\LDER^q\km - \lder^p\bigl(\mu\kmt -
 \theta(\mu)\ko\bigr)\otimes\LDER^q\bigl((\km -
 \ko)/\bar{\mu}\bigr)\\
 &= \lder^p\kmt\otimes\LDER^q\km -
 \lder^p\bigl(\mu\kmt\bigr)\otimes\LDER^q\left
 (\km/\bar{\mu}\right)
 + Z_{pq},
 \end{split}
 \end{equation}
 where $Z_{pq}$ is an operator of the form $\psi_{pq}\otimes\ko +
 \ko\otimes\chi_{pq}$. Using the identity $(zf)^{(p)} =
 pf^{(p-1)} + z f^{(p)}$, we get
 $$
 \lder^p\bigl(\mu\kmt\bigr) = p\lder^{p-1}\kmt +
 \mu\lder^p\kmt, \qquad  \LDER^q\km =
  \LDER^q\bigl(\bar\mu(\km/\bar\mu)\bigr) = q
 \LDER^{q-1}\bigl(\km/\bar\mu\bigr) +
 \bar\mu\LDER^q\bigl(\km/\bar\mu\bigr).
 $$
 Substituting this into \eqref{eq-20}, we obtain
 \begin{equation}\label{eq-22}
 T_{pq} - \Sth T_{pq}\Sth^* = q
 \left(\lder^p\kmt\otimes\LDER^{q-1}
 \bigl(\km/\bar\mu\bigr)\right)
 - p \left(\lder^{p-1}\kmt\otimes\LDER^q\bigl
 (\km/\bar\mu\bigr)\right)
 + Z_{pq}.
 \end{equation}
 For the operators $T_{00}$, $T_{01}$ and $T_{10}$ we have
 \begin{equation}\label{eq-50}
 \begin{split}
 &T_{00} - \Sth T_{00}\Sth^* = Z_{00};\\
 &T_{01} - \Sth T_{01}\Sth^* = \kmt\otimes \left(\km/\bar\mu\right) + Z_{01};\\
 &T_{10} - \Sth T_{10}\Sth^* = -\kmt\otimes\left(\km/\bar\mu\right) + Z_{10}.
 \end{split}
 \end{equation}
 Since $A$ is a truncated Toeplitz operator, it satisfies \eqref{s}
 with some $\psi,\chi \in \Kth$. Combining \eqref{s}, \eqref{eq-22} and
 \eqref{eq-50}, we obtain
 \begin{equation}\label{eq-23}
 \Psi\otimes\ko + \ko\otimes\Phi = \sum_{0 \le p,q\le n}
 \bigr((q+1)a_{p, q+1} - (p+1)a_{p+1,q}\bigr)
 \cdot\left(\lder^{p}\kmt\otimes\LDER^{q}
 \bigl(\km/\bar\mu\bigr)\right),
 \end{equation}
 where $\Psi,\Phi \in \Kth$ and $a_{n+1,q} = a_{p,n+1} = 0$ for all $0
 \le p,q \le n$. It follows from Proposition \ref{prp14} that
 \begin{equation}\label{eq-24}
 (q+1)a_{p, q+1} - (p+1)a_{p+1,q} = 0, \quad 0 \le p,q \le n.
 \end{equation}
 In the above formula there is no restriction on $a_{0,0}$, which agrees well with the fact that $ \kmt\otimes\km \in \Tau$. For each $1 \le s \le n$, from \eqref{eq-24} we get  the following system:
 \begin{equation}
 \left.
 \begin{aligned}
 sa_{0,s} - a_{1, s-1}&=0\\
 (s-1)a_{1,s-1}- 2a_{2, s-2}& = 0\\
 \ldots\\
 2a_{s-2,1}- (s-1)a_{s-1,1}& = 0\\
 a_{s-1,1}- sa_{s,0}& = 0
 \end{aligned}
 \right\}
 \end{equation}
 Solving this system, we obtain
$$
 a_{t, s-t} = a_{0,s}C_{s-t}^{t}, \quad C_{s-t}^{t} = \frac{s!}{t!(s-t)!}, \quad 0
 \le t \le s.
$$
 It follows from \eqref{eq-24} that $na_{n,n} - (n+1)a_{n+1, n-1} = 0$ and thus $a_{n,n} = 0$. By induction, we have $a_{p,q} = 0$ for all indexes $p,q$ such that  $p + q > n$.
 Now we get the required representation from formulas \eqref{eq-31} and
  \eqref{eq-19}:
\begin{equation}\label{eq-80}
 A = \sum_{s =0}^{n} \sum_{t=0}^{s}
 a_{t,s-t}\left(\lder^t\kmt\otimes\LDER^{s-t}\km\right) =
 \sum_{s =0}^{n} a_{0,s} \cdot \der^s [\kmt\otimes\km].
 \end{equation}
In the case $\mu =0$, put $T_{pq} =  \lder^p\kot\otimes\LDER^q\ko$. It follows from Proposition \ref{lem-nder} that
 \begin{equation*}
 \begin{split}
 &T_{pq} - \Sth T_{pq}\Sth^* = \lder^p\kot\otimes\LDER^q\ko -
 \frac{p}{q+1}\lder^{p-1}\kot\otimes\LDER^{q+1}\ko + Z_{pq}, \quad p \ge 1;\\
 &T_{pq} - \Sth T_{pq}\Sth^* = \lder^p\kot\otimes\LDER^q\ko + Z_{pq},
 \quad p = 0.
 \end{split}
 \end{equation*}
Proceeding as in the case $\mu \neq 0$, we obtain the system
 $$a_{p,q} - \frac{p+1}{q}a_{p+1,q-1} = 0, \qquad 0 \le p \le n, \quad 1 \le q \le n+1,$$
where $a_{n+1,q} = 0$ for all $q$ and $a_{p,-1}=0$ for all $p$. This system has the same solution as the system in \eqref{eq-24}. Hence we have representation \eqref{eq-80} in the case $\mu = 0$ as well. \qed

\enddocument